\documentclass{article}

\usepackage{
 amsmath,
 amsxtra,
 amsthm,
 amssymb,
 mathrsfs,
 stmaryrd
  }
\usepackage[all]{xy}

\newtheorem{theorem}{Theorem}[section]
\newtheorem{lemma}[theorem]{Lemma}

\newtheorem{proposition}[theorem]{Proposition}
\newtheorem{corollary}[theorem]{Corollary}
\newtheorem*{thm2}{Theorem}

\theoremstyle{definition}
\newtheorem{defn}[theorem]{Definition}
\newtheorem{remark}[theorem]{Remark}

\newcommand{\bd}{\begin{defn}}
\newcommand{\ed}{\end{defn}}
\newcommand{\bl}{\begin{lemma}}
\newcommand{\el}{\end{lemma}}
\newcommand{\bp}{\begin{proposition}}
\newcommand{\ep}{\end{proposition}}
\newcommand{\bt}{\begin{theorem}}
\newcommand{\et}{\end{theorem}}
\newcommand{\bc}{\begin{corollary}}
\newcommand{\ec}{\end{corollary}}
\newcommand{\br}{\begin{remark}}
\newcommand{\er}{\end{remark}}
\newcommand{\ba}{\begin{array}}
\newcommand{\ea}{\end{array}}
\newcommand{\bpf}{\begin{proof}}
\newcommand{\epf}{\end{proof}}

\newcommand{\Qp}{\mathbb{Q}_p}
\newcommand{\Z}{\mathbb{Z}}
\newcommand{\Zp}{\mathbb{Z}_p}
\newcommand{\Op}{\mathcal{O}}
\newcommand{\Ga}{\Gamma}
\newcommand{\ga}{\gamma}

\DeclareMathOperator{\Gal}{Gal}
\DeclareMathOperator{\Hom}{Hom}

\DeclareMathOperator{\rank}{rank}
\DeclareMathOperator{\Ext}{Ext}

\newcommand{\Iw}{\mathrm{Iw}}

\newcommand{\Tr}{\mathrm{Tr}}

\newcommand{\lra}{\longrightarrow}

\newcommand{\ot}{\otimes}

\newcommand{\ps}[1]{\llbracket #1 \rrbracket}

\newcommand{\ilim}{\displaystyle \mathop{\varinjlim}\limits}
\newcommand{\plim}{\displaystyle \mathop{\varprojlim}\limits}

\numberwithin{equation}{section}

\begin{document}
\title{Comparing direct limit and inverse limit of even $K$-groups in non-commutative $p$-adic Lie extensions}
 \author{
  Meng Fai Lim\footnote{School of Mathematics and Statistics,
Central China Normal University, Wuhan, 430079, P.\ R.\ China.
 E-mail: \texttt{limmf@ccnu.edu.cn}} }
\date{}
\maketitle

\begin{abstract} \footnotesize
\noindent In a previous paper of the author, we establish a duality for the direct limit and the inverse limit of higher even $K$-groups over a $\Zp^d$-extension. In this paper, we shall establish such a duality over certain non-commutative $p$-adic Lie extensions.


\end{abstract}

\section{Introduction} \label{intro}

In this paper, $p$ will always denote an odd prime. We shall also fix an integer $i\geq 2$. Let $F$ be a number field and $F_\infty$ a $\Zp^d$-extension of $F$.  For each finite intermediate extension $L$ of $F_\infty/F$, the works of Quillen \cite{Qui73b} and Borel \cite{Bo} tell us that the even $K$-groups $K_{2i-2}(\Op_L)$ are finite. Now, for two finite subextensions $L\subseteq L'$, the inclusion $\Op_L\lra \Op_{L'}$ induces a map 
\[\jmath_{L/L'}: K_{2i-2}(\Op_L)[p^\infty]\lra K_{2i-2}(\Op_{L'})[p^\infty]\]
 by functoriality. (Note that this induced map may not be injective in general.) In the other direction, there is the norm map (also called the transfer map) 
 \[\Tr_{L'/L}: K_{2i-2}(\Op_{L'})[p^\infty]\lra K_{2i-2}(\Op_{L})[p^\infty].\]
  We then consider the following direct limit and inverse limit
\[  \ilim_L K_{2i-2}(\Op_L)[p^\infty] \quad \mbox{and}\quad \plim_L K_{2i-2}(\Op_L)[p^\infty],\]
whose transition maps are given by the maps $\jmath_{L/L'}$ and $\Tr_{L'/L}$ respectively. These limit modules come equipped with natural $\Zp\ps{G}$-module structures, where $G=\Gal(F_\infty/F)\cong \Zp^d$. In \cite[Theorem 3.1]{LimKlimit}, the author proved the following.

\begin{thm2}[\cite{LimKlimit}]
Retain the notation as above. Then there is a pseudo-isomorphism
\[  \left(\plim_L K_{2i-2}(\Op_L)[p^\infty]\right)^\iota \sim\left(\ilim_L K_{2i-2}(\Op_L)[p^\infty]\right)^\vee\]
of $\Zp\ps{G}$-modules.
\end{thm2}

Here, for a given $\Zp\ps{G}$-module $M$, the module $M^\iota$ is defined to be the same underlying $\Zp$-module $M$ but with a $\Ga$-action given by
\[ \ga \cdot_{\iota} x = \ga^{-1}x, \quad \gamma \in\Ga, x\in M,\]  
and $(-)^\vee$ is the Pontryagin dual.
Since $K_0(\Op_L)[p^\infty]$ can be identified with the Sylow $p$-subgroup of the class group of $\Op$, the above result may be thought as a generalization of the previous results on class groups (see the works of Nekov\'a\v{r} \cite{Ne}, Vauclair \cite{Vau} and, more recently, that of Lai and Tan \cite{LT}) to the higher even $K$-groups. 

The goal of this paper is to extend the preceding theorem to certain non-commutative $p$-adic Lie
extensions. As before, $F$ will always denote a number field. We shall, once and for all, fix a finite set $S$ of primes of $F$ which is always assumed to contain those above $p$ and the infinite primes. Let $S'$ denote the set of finite primes of $S$. Write $F_S$ for the maximal algebraic extension of $F$ unramified outside $S$. For each finite extension $L$ of $F$ contained in $F_S$, denote by $S(L)$ (resp., $S'(L)$) the set of primes of $L$ above $S$ (resp., above $S'$). For ease of notation, we shall write $\Op_{L,S}$ (instead of $\Op_{L,S'}$) for the ring of $S'(L)$-integers of $L$. This notational convention will be in force throughout the remainder of the paper without further mention. 

For two finite subextensions $L\subseteq L'$, we have a map $\jmath_{L/L'}: K_{2i-2}(\Op_{L,S})[p^\infty]\lra K_{2i-2}(\Op_{L',S})[p^\infty]$ induced by the inclusion $\Op_{L,S}\lra \Op_{L',S}$ and a norm map $\Tr_{L'/L}: K_{2i-2}(\Op_{L',S})[p^\infty]\lra K_{2i-2}(\Op_{L,S})[p^\infty]$. Similarly, we consider the following direct limit and inverse limit
\[  \ilim_L K_{2i-2}(\Op_{L,S})[p^\infty] \quad \mbox{and}\quad \plim_L K_{2i-2}(\Op_{L,S})[p^\infty],\]
whose transition maps are given by the maps $\jmath_{L/L'}$ and $\Tr_{L'/L}$ respectively.

Our main result is as follows.

\bt \label{Klimitthm}
Let $i\geq 2$ and let $F_\infty$ be a $p$-adic Lie extension of $F$ which is contained in $F_S$. Assume that $G=\Gal(F_\infty/F)$ has no $p$-torsion. Suppose that at least one of the following statements holds.
\begin{enumerate}
  \item[$(a)$] The integer $i$ is even and $F_\infty$ is a totally real field.
  \item[$(b)$] $F_\infty = F(\mu_{p^\infty}, a_1^{p^{-\infty}},..., a_r^{p^{-\infty}})$.
\end{enumerate}
Then there is an isomorphism
\[  \left(\ilim_L K_{2i-2}(\Op_{L,S})[p^\infty]\right)^\vee \cong \Ext^1_{\Zp\ps{G}}\left(\plim_L K_{2i-2}(\Op_{L,S})[p^\infty], \Zp\ps{G}\right)\]
of $\Zp\ps{G}$-modules.
\et

To see how our theorem compares with that in \cite{LimKlimit}, we first return to the situation when $G=\Gal(F_\infty/F)\cong \Zp^d$. In this context, it follows from \cite[Proposition 8]{PR} that $\Ext^1_{\Zp\ps{G}}(M, \Zp\ps{G})$ is pseudo-isomorphic to $M^\iota$. Furthermore, \cite[Theorem 1]{Iw73} tells us that $F_\infty \subseteq F_{S_{p\infty}}$, where $S_{p\infty}$ is the set of primes consisting precisely of those above $p$ and the infinite primes.
On the other hand, by an application of the localization sequence of Soul\'e (see \cite[Proposition 2.2]{Kol} or \cite[Section III.3]{Sou}), one can show that $K_{2i-2}(\Op_{L,S_p})[p^\infty] = K_{2i-2}(\Op_{L})[p^\infty]$ for every $L\subseteq F_\infty$. In view of these observations, the main result in \cite{LimKlimit} can be restated as the following isomorphism
\[  \left(\ilim_L K_{2i-2}(\Op_{L,S_p})[p^\infty]\right)^\vee \cong \Ext^1_{\Zp\ps{G}}\left(\plim_L K_{2i-2}(\Op_{L,S_p})[p^\infty], \Zp\ps{G}\right).\]
This is the form that we will prove in our noncommutative settings.

We end the introductory section giving an outline of the paper. Section \ref{Iwasawa modules} reviews certain preliminaries on Iwasawa algebras and modules. Section
\ref{K and Galois} is where we recall the relation between
the Sylow $p$-subgroups of the even $K$-groups and Galois cohomology
groups. We also introduce the Iwasawa cohomology group here and end the section with a proof of our main theorem under hypothesis (a). Section \ref{false-Tate section}, which is the final section of the paper, will prove the main theorem under hypothesis (b).

\section{Iwasawa modules} \label{Iwasawa modules}

As before, we let $p$ denote a fixed odd prime. The Iwasawa algebra of a compact $p$-adic Lie group $G$ over $\Zp$ is defined by 
\[ \Zp\ps{G} = \plim_U \Zp[G/U], \]
where $U$ runs over the open normal subgroups of $G$ and the transition maps in the inverse
limit are given by the canonical projection maps.

In the event that the group $G$ is pro-$p$ with no $p$-torsion, it is known that $\Zp\ps{G}$ is
a Noetherian Auslander regular ring (cf.\ \cite[Theorem 3.26]{V02}) with no zero divisors (cf.\
\cite{Neu}). As a consequence, it admits a skew field $Q(G)$ which is flat
over $\Zp\ps{G}$ (see \cite[Chapters 6 and 10]{GW} or \cite[Chapter
4, \S 9 and \S 10]{Lam}). This in turn enables us to define the $\Zp\ps{G}$-rank of a finitely generated $\Zp\ps{G}$-module $M$, which is given by
$$ \rank_{\Zp\ps{G}}(M)  = \dim_{Q(G)} (Q(G)\ot_{\Zp\ps{G}}M). $$
The module $M$ is then said to be 
torsion over $\Zp\ps{G}$ if $\rank_{\Zp\ps{G}} (M) = 0$. We also recall a useful equivalent formulation of a torsion $\Zp\ps{G}$-module: $\Hom_{\Zp\ps{G}}(M,\Zp\ps{G})=0$ (cf.\ \cite[Remark 3.7]{V02}). If the torsion $\Zp\ps{G}$-module $M$ also satisfies $\Ext^1_{\Zp\ps{G}}(M,\Zp\ps{G})=0$, we shall say that $M$ is a pseudo-null $\Zp\ps{G}$-module.

We now extend the notion of torsion modules and pseudo-null modules to the case when $G$ is a general compact $p$-adic Lie group. A well-known theorem of Lazard asserts that $G$ contains an open normal subgroup $G_0$ which is pro-$p$ with no $p$-torsion (cf.\ \cite[Theorem 8.32]{DSMS}). By \cite[Proposition 5.4.17]{NSW}, we have
\[\Ext^i_{\Zp\ps{G}}(M,\Zp\ps{G}) \cong \Ext^i_{\Zp\ps{G_0}}(M,\Zp\ps{G_0})\]
for every finitely generated $\Zp\ps{G}$-module $M$. In view of this identification, we shall say that $M$ is a torsion $\Zp\ps{G}$-module (resp., psuedo-null $\Zp\ps{G}$-module) if $\Hom_{\Zp\ps{G}}(M,\Zp\ps{G})=0$ (resp., $\Ext^i_{\Zp\ps{G}}(M,\Zp\ps{G})=0$ for $i=0,1$). Equivalently, this is saying that $M$ is a torsion $\Zp\ps{G}$-module (resp., pseudo-null $\Zp\ps{G}$-module), whenever $M$ is a torsion $\Zp\ps{G_0}$-module (resp., pseudo-null $\Zp\ps{G_0}$-module) in the sense of the preceding paragraph. (Also, compare with \cite[Discussion after Definition 2.6]{V02}.)

\bd
Let $M$ be a $\Zp\ps{G}$-module. For an open subgroup $U$ of $G$, we write $M_U$ for the largest quotient
of $M$ on which $U$ acts trivially. We shall say that the $\Zp\ps{G}$-module $M$ is \textit{systematically coinvariant-finite} if $M_U$ is finite for every open subgroup $U$ of $G$. Note that $\{M_U\}_U$ forms a direct system of finite modules with transition maps given by the norm maps \[ N_{U/V}: M_U\lra M_V\]
for $V\subseteq U$. In particular, $\ilim_U M_U$ is a discrete $\Zp\ps{G}$-module if $M$ is systematically coinvariant-finite.
\ed

\bl \label{finite fg tor}
A systematically coinvariant-finite $\Zp\ps{G}$-module $M$ is finitely generated torsion over $\Zp\ps{G}$.
\el

\bpf
Replacing $G$ by a smaller subgroup, we may assume that $G$ is pro-$p$ with no $p$-torsion. Since $M_G$ is finite, Nakayama lemma tells us that $M$ is finitely generated over $\Zp\ps{G}$. It then follows from a formula of Harris \cite[Theorem 1.10]{Har} that
\[\rank_{\Zp}(M_{G_n}) = \rank_{\Zp\ps{G}}(M) p^{dn} +O(p^{d(n-1)}),\]
where $d$ is the dimension of $G$ (in the sense of \cite[Theorem 8.36]{DSMS}) and $\{G_n\}$ is a certain prescribed decreasing sequence of open subgroups with $\cap_{n\geq 1} G_n =\{1\}$ (see \cite[P. 105]{Har} for its precise definition).
In view of the systematically coinvariant-finiteness assumption, we have $\rank_{\Zp}(M_{G_n}) =0$ for every $n$. This in turn forces $\rank_{\Zp\ps{G}}(M)=0$ in the above equality which proves the lemma.
\epf

The next short exact sequence will be a key component for our subsequent discussion.

\bp \label{alg main}
For every systematically coinvariant-finite $\Zp\ps{G}$-module $M$, we have the following short exact sequence
\[0 \lra \Big(\ilim_U M_U\Big)^\vee \lra \Ext^1_{\Zp\ps{G}}(M,\Zp\ps{G})
           \lra \Big(\ilim_U H_1(U,M)\ot\Qp/\Zp\Big)^\vee \lra 0,\]
where the transition maps in the limit is given by norm maps.
\ep

\bpf By a result of Jannsen (see \cite[Theorem 2.1]{Jannsen89} or \cite[Theorem 5.4.13]{NSW}), there is a short exact sequence
\[0 \lra \Big(\ilim_U(M_U[p^\infty])\Big)^\vee \lra \Ext^1_{\Zp\ps{G}}(M,\Zp\ps{G})
           \lra \Big(\ilim_U H_1(U,M)\ot\Qp/\Zp\Big)^\vee \lra 0. \]
           On the other hand, the systematically coinvariant-finite property implies that the leftmost term in the above exact sequence is $\Big( \ilim_U M_U\Big)^\vee$. 
\epf

\section{$K$-groups and Galois cohomology} \label{K and Galois}

Retain the notation from Section \ref{intro}. We now recall the following deep calculations of Qullien \cite{Qui73a} and Borel \cite{Bo}.

\bt \label{Qullien-Borel}
Let $i\geq 2$. Then the group $K_{2i-2}(\Op_{F,S})$ is finite and 
\[ \rank_\Z\big( K_{2i-1}(\Op_{F,S})\big)= \left\{
  \begin{array}{ll}
    r_1(F)+r_2(F), & \hbox{if $i$ is odd;} \\
    r_2(F), & \hbox{if $i$ is even.}
  \end{array}
\right. \]
Here $r_1(F)$ (resp., $r_2(F)$) is the number of real embeddings (resp., number of
pairs of complex embeddings) of $L$.
\et

Denote by $\mu_{p^\infty}$ the group of all $p$-power root of unity. The natural action of $G_S(F)$ on $\mu_{p^\infty}$ induces the cyclotomic character
\[\chi: G_S(F) \lra \mathrm{Aut}(\mu_{p^\infty}) \cong \Zp^{\times}.\]
For a $\Zp[G_S(F)]$-module $X$, we denote by $X(i)$ the $i$-fold Tate twist of $X$. More precisely, $X(i)$ is the $G_S(F)$-module which is $X$ as a $\Zp$-module but with a $G_S(F)$-action given by
\[ \sigma\cdot x = \chi(\sigma)^i\sigma x,  \]
where the action on the right is the original action of $G_S(F)$ on $X$.
Plainly, we have $X(0)=X$ and $\mu_{p^{\infty}} \cong \Qp/\Zp(1)$.
In \cite{Sou}, Soul\'e connected the higher $K$-groups with (continuous) Galois cohomology groups via the $p$-adic Chern class maps
\[ \mathrm{ch}_{i,k}^{(p)}: K_{2i-k}(\Op_{F,S})\ot \Zp \lra H^k(G_{S}(F), \Zp(i))\]
for $i\geq 2$ and $k =1,2$. (For the precise definition of these maps, we refer readers to loc.\ cit.) By the deep work of Rost and Voevodsky \cite{Vo} (also see \cite{Wei09}), these $p$-adic Chern class maps are known to be isomorphisms. As a consequence of these, we have the following identification between the Sylow $p$-subgroup of $K$-groups and Galois cohomology.

\bp \label{K2 = H2}
For $i\geq 2$, one has the following isomorphisms 
\[ K_{2i-2}(\Op_{F,S})[p^\infty] \cong H^2\big(G_{S}(F), \Zp(i)\big)\]
of finite abelian groups. Furthermore, if $F$ is totally real and $i$ is even, then one has the following isomorphisms 
\[ K_{2i-1}(\Op_{F,S})[p^\infty] \cong H^1\big(G_{S}(F), \Zp(i)\big)\]
of finite abelian groups.
\ep

We now introduce the Iwasawa cohomology groups which will be an important object for our study. For every extension $L$ of $F$ contained in $F_S$, write $G_S(L)$ for the Galois group $\Gal(F_S/L)$.
For every extension $\mathcal{L}$ of $F$ contained in $F_S$, the Iwasawa cohomology groups are defined by
 \[H^k_{\Iw,S}\big(\mathcal{L}/F, \Zp(i)\big):= \plim_L H^k\big(G_S(L), \Zp(i)\big), \]
 where $L$ runs through all finite extensions of $F$ contained in $\mathcal{L}$ and the transition maps are given by the corestriction maps. For ease of notation,
we will drop the `$S$' in the Iwasawa cohomology groups.

\br \label{finite Iw}
Note that if $\mathcal{L}/F$ is a finite extension, then one has
 $H^k_{\Iw,S}\big(\mathcal{L}/F, \Zp(i)\big)= H^k\big(G_S(\mathcal{L}), \Zp(i)\big)$.
\er 

Let $F_\infty$ be a $p$-adic Lie extension of $F$ contained in $F_S$. In other words, $F_\infty$ is a Galois extension of $F$, whose Galois group $G=\Gal(F_\infty/F)$ is a compact $p$-adic Lie group. The groups $H^k_{\Iw}\big(F_\infty/F, \Zp(i)\big)$ come naturally equipped with $\Zp\ps{G}$-module structures and can be shown to be finitely generated over $\Zp\ps{G}$ (see \cite[Proposition 4.1.3]{LimSh}). 
The following version of Tate's descent spectral sequence for Iwasawa cohomology will be useful in our subsequent discussion.

\bp \label{gal descent}
Let $F_\infty$ be a $p$-adic Lie extension of $F$ contained in $F_S$. Suppose that $i\geq 2$. 
Let $U$ be a closed normal subgroup of $G=\Gal(F_\infty/F)$ and write $\mathcal{L}$ for the fixed field of $U$. Then we have a homological spectral sequence
 $$ H_r\big(U, H^{2-s}_{\Iw}(F_\infty/F, \Zp(i))\big)\Longrightarrow H^{2-r-s}_{\Iw}\big(\mathcal{L}/F, \Zp(i)\big) $$
and an isomorphism
\[ H^2_{\Iw}\big(F_\infty/F, \Zp(i)\big)_U \cong H^2_{\Iw}\big(\mathcal{L}/F, \Zp(i)\big).\]
Furthermore, if $U$ is an open subgroup of $G$, then $\mathcal{L}$ is a finite extension of $F$ and we have
\[ H^2_{\Iw}\big(F_\infty/F, \Zp(i)\big)_U \cong H^2\big(G_S(\mathcal{L}), \Zp(i)\big).\]
\ep

\bpf
When $G$ is a commutative group, the spectral sequence is a result of Nekov\'a\v{r} \cite[Proposition 4.2.3]{Ne}. For a noncommutative group $G$, this follows from \cite[Proposition 1.6.5]{FK} or \cite[Theorem 3.1.8]{LimSh}. The isomorphisms of the proposition follow from reading off the initial term of the spectral sequence. Finally, in the event that $U$ is an open subgroup of $G$, then $\mathcal{L}$ is a finite extension of $F$, and so $H^2_\Iw(L/F,\Zp(i))$ identifies with $H^2(G_S(L), \Zp(i))$ as noted in Remark \ref{finite Iw}.
\epf

We now record the following which was also proved in \cite[Proposition 4.1.1]{LimKgroups}.

\bp \label{tor}
Let $i\geq 2$. For every $p$-adic Lie extension $F_\infty$ of $F$, the module $H^2_{\Iw}\big(F_\infty/F, \Zp(i)\big)$ is torsion over $\Zp\ps{\Gal(F_\infty/F)}$.
\ep

\bpf
We shall give a different proof to that in \cite[Proposition 4.1.1]{LimKgroups}. From Proposition \ref{gal descent}, we see that for every open subgroup $U$ of $G=\Gal(F_\infty/F)$,  \[ H^2_{\Iw}\big(F_\infty/F, \Zp(i)\big)_U \cong H^2\big(G_S(\mathcal{L}), \Zp(i)\big) \] where the latter is finite in view of Proposition \ref{K2 = H2}. Therefore,  $H^2_{\Iw}\big(F_\infty/F, \Zp(i)\big)$ is a systematically coinvariant-finite $\Zp\ps{G}$-module, and so is necessarily torsion over $\Zp\ps{G}$ by Lemma \ref{finite fg tor}.
\epf

We now apply Proposition \ref{Jannsenses} to $H^2_{\Iw}\big(F_\infty/F, \Zp(i)\big)$. We first recall that from the discussion in \cite[Chap. III]{Sou}, the isomorphism in Proposition \ref{K2 = H2} fits into the following commutative diagram
\begin{equation} \label{N cor}\entrymodifiers={!! <0pt, .8ex>+} \SelectTips{eu}{}
\xymatrix{
      K_{2i-2}(\Op_{L,S})[p^\infty] \ar[r]^{\mathrm{ch}_i^L} \ar[d]^{\Tr_{L/F}} &  H^2\left(G_S(L), \Zp(i)\right) \ar[d]^{\mathrm{cor}} \\
      K_{2i-2}(\Op_{F,S})[p^\infty] \ar[r]^{\mathrm{ch}_i^F} &  H^2\left(G_S(F), \Zp(i)\right)}
      \end{equation}
      \begin{equation} \label{j tr}\entrymodifiers={!! <0pt, .8ex>+} \SelectTips{eu}{}
\xymatrix{
      K_{2i-2}(\Op_{F,S})[p^\infty] \ar[r]^{\mathrm{ch}_i^L} \ar[d]^{j_{L/F}} &  H^2\left(G_S(F), \Zp(i)\right) \ar[d]^{\mathrm{res}} \\
      K_{2i-2}(\Op_{L,S})[p^\infty] \ar[r]^{\mathrm{ch}_i^F} &  H^2\left(G_S(L), \Zp(i)\right)}
      \end{equation}
From the diagram (\ref{N cor}), we obtain
\[ H^2_{\Iw}\big(F_\infty/F, \Zp(i)\big) \cong  \plim_{L} K_{2i-2}(\Op_{L,S})[p^\infty]. \]

On the other hand, by a similar argument to that in \cite[Lemma 3.6]{LimKlimit}, one has the following commutative diagram
\[ \entrymodifiers={!! <0pt, .8ex>+} \SelectTips{eu}{}\xymatrix{
      H^2_{\Iw}\big(F_\infty/F, \Zp(i)\big)_U \ar[r]_{\sim}\ar[d]_{N_{U/V}} &  H^2\big(G_S(F_\infty^U), \Zp(i)\big) \ar[d]_{\mathrm{res}}\\
      H^2_{\Iw}\big(F_\infty/F, \Zp(i)\big)_V  \ar[r]_{\sim} &  H^2\big(G_S(F_\infty^V), \Zp(i)\big)
     }\]
     for every open normal subgroups $V\subseteq U$ of $G$, where $F_\infty^U$ (resp., $F_\infty^V$) denotes the fixed field of $F_\infty$ under $U$ (resp., under $V$). Upon combining this with (\ref{j tr}), we have
\[  \ilim_{U, N_{U/V}}  H^2_{\Iw}\big(F_\infty/F, \Zp(i)\big)_U \cong  \ilim_{L} K_{2i-2}(\Op_{L,S})[p^\infty].\]

Putting all the above observations together into Proposition \ref{Jannsenses}, we obtain a short exact sequence
\begin{multline}\label{Jannsenses}
   0 \lra \Big(\ilim_L K_{2i-2}(\Op_{L,S})[p^\infty]\Big)^\vee \lra \Ext^1_{\Zp\ps{G}}\Big(\plim_{L} K_{2i-2}(\Op_{L,S})[p^\infty],\Zp\ps{G}\Big) \\
 \lra \Big(\ilim_U H_1\big(U,H^2_{\Iw}\big(F_\infty/F, \Zp(i)\big)\big)\ot\Qp/\Zp\Big)^\vee \lra 0.
\end{multline}

Therefore, it remains to show that $H_1\big(U, H^2_{\Iw}\big(F_\infty/F, \Zp(i)\big)\big)$ is finite for every open subgroup $U$ of $G$. Without loss of generality, upon relabelling, it is enough to show that $H_1\big(G, H^2_{\Iw}\big(F_\infty/F, \Zp(i)\big)\big)$ is finite. 

We now give a proof of this in the totally real context.

\bpf[Proof of Theorem \ref{Klimitthm} under assumption (a)] From the low degree terms of the homological spectral sequence
 $$ H_r\big(G, H^{2-s}_{\Iw}(F_\infty/F, \Zp(i))\big)\Longrightarrow H^{2-r-s}\big(G_S(F), \Zp(i)\big), $$
 we have a surjection
\[   H^{1}\big(G_S(F), \Zp(i)\big) \lra H_1\big(G, H^2_{\Iw}(F_\infty/F, \Zp(i))\big) \lra 0. \]
Since $F$ is totally real and $i$ is even, it follows from Theorem \ref{Qullien-Borel} and Proposition \ref{K2 = H2} that the group $H^{1}\big(G_S(F), \Zp(i)\big)$ is finite. Hence it follows that 
$H_1\big(G, H^2_{\Iw}\big(F_\infty/F, \Zp(i)\big)\big)$ is also finite which is what we want to show.
\epf

\section{Multi-false Tate extension} \label{false-Tate section}

The goal of this final section of the paper is to prove Theorem \ref{Klimitthm} under assumption (b). Set $F_\infty = F(\mu_{p^\infty}, a_1^{p^{-\infty}},..., a_r^{p^{-\infty}})$, where $a_1,..., a_r \in F^\times$ whose image in $F^\times/(F^\times)^p$ are linearly independent over $\Z/p\Z$. (This is usually coined as the multi-false Tate extension.) We then write $G=\Gal(F_\infty/F)$ and $H= \Gal(F_\infty/ F(\mu_{p^\infty}))$, where Kummer theory tells us that $H\cong \Zp^r$. As seen in the discussion in the previous section, it remains to show that $H_1\big(G, H^2_{\Iw}\big(F_\infty/F, \Zp(i)\big)\big)$ is finite. Now, set $G_0 = \Gal(F_\infty/F(\mu_p))$. Clearly, $G_0$ is an open normal subgroup of $G$ with index coprime to $p$. Therefore, one has an identification
\[  H_1\big(G_0, H^2_{\Iw}\big(F_\infty/F, \Zp(i)\big)\big)_{G/G_0}\cong  H_1\big(G, H^2_{\Iw}\big(F_\infty/F, \Zp(i)\big)\big). \]
In view of this isomorphism, for the verification of the finiteness of $H_1\big(G, H^2_{\Iw}\big(F_\infty/F, \Zp(i)\big)\big)$, one is reduced to showing that $H_1\big(G_0, H^2_{\Iw}\big(F_\infty/F, \Zp(i)\big)\big)$ is finite. Therefore, replacing $F$ by $F(\mu_p)$, we may assume that the group $\Gamma := G/H = \Gal(F(\mu_{p^\infty})/F)$ is isomorphic to $\Zp$ which we will do for the remainder of the section.

We shall require a further technical result. Before stating this said result, we make a remark on notation. 
For a given character $\kappa : G/H \lra
\Zp^{\times}$ and a $\Zp\ps{G/H}$-module $T$, we let
$T(\kappa)$ denote the $\Zp$-module $T$ with the new commuting
$G/H$-action given by the twist of the original by $\kappa$. With this notation in hand, we can state the required technical result.

\bp \label{Kato} Let $G$ be a compact $p$-adic Lie group and $H$
a closed normal subgroup of $G$. Assume that we are given a
finite family of closed normal subgroups $H_j$ $(0\leq j\leq r)$ of
$G$ such that $1=H_0\subseteq H_1 \subseteq \cdots\subseteq H_r =H$ with
$H_j/H_{j-1}\cong \Zp$ for $1\leq j\leq r$, and such that the action
of $G$ on $H_j/H_{j-1}$ by inner automorphism is given by a
homomorphism $\chi_j:G/H\lra \Zp^{\times}$.

For a given finitely generated $\Zp\ps{G}$-module $M$, there exists a
finite family $(S_{k})_{1\leq k\leq t}$ of
$\Zp\ps{G/H}$-submodules of $H_1(H,M)$ satisfying all of the following
properties. \vspace{-0.05in}
\begin{enumerate}
\item[$(i)$] $0 = S_0\subseteq S_1\subseteq \cdots\subseteq S_t = H_1(H,M)$.

\vspace{-0.05in}
\item[$(ii)$] For each $k$ $(1\leq k\leq t)$, there is a
$\Zp\ps{G/H}$-subquotient $T$ of $M_H$ and a family
$(s(k))_{1\leq k\leq r}$ of nonnegative
integers such that $|\{k |s(k)>0\}|\geq 1$ and such that
$S_k/S_{k-1}$ is isomorphic to the twist $T(\prod_{1\leq i\leq
r}\chi_k^{s(k)})$ of $T$.\end{enumerate} \ep

\bpf This is a special case of a result of Kato \cite[Proposition 4.2]{Ka2}, and we refer the readers to loc. cit for the proof. \epf

We return to our multi-false Tate extension situation. Kummer theory tells us that the action of $\Gamma$ on $H$ by inner automorphism coincides with the cyclotomic character. Therefore, in this context, each $S_k/S_{k-1}$ is isomorphic to $T(n_k)$ for some positive integer $n_k$.

On the other hand, from the spectral sequence
\[ H_r\big(\Gamma, H_s\big(H,  H^2_{\Iw}\big(F_\infty/F, \Zp(i)\big)\big)\big) \Longrightarrow H_{r+s}\big(G, H^2_{\Iw}\big(F_\infty/F, \Zp(i)\big)\big),  \]
we obtain a short exact sequence
\[ 0\lra H_1\big(H,  H^2_{\Iw}\big(F_\infty/F, \Zp(i)\big)\big)_{\Ga}\lra  H_{1}\big(G, H^2_{\Iw}\big(F_\infty/F, \Zp(i)\big)\big)\lra H_1\big(\Gamma, H^2_{\Iw}\big(F_\infty/F, \Zp(i)\big)_H\big)\lra 0.  \]

Therefore, we are now reduced to showing the leftmost and rightmost terms of the short exact sequence is finite. By Proposition \ref{gal descent}, we have 
\[ H^2_{\Iw}\big(F_\infty/F, \Zp(i)\big)_H\cong H^2_{\Iw}\big(F(\mu_{p^\infty})/F, \Zp(i)\big)\]
which is a torsion $\Zp\ps{\Gamma}$-module. Now by Proposition \ref{gal descent} again, we have
\[ H^2_{\Iw}\big(F(\mu_{p^\infty})/F, \Zp(i)\big)_\Gamma\cong H^2\big(G_S(F), \Zp(i)\big)\]
which is finite. In view of this latter finiteness observation, it follows from the structure theory of finitely generated $\Zp\ps{\Gamma}$-module (for instances, see \cite[Proposition A.1.7]{CS}) that $H_1\big(\Gamma, H^2_{\Iw}\big(F(\mu_{p^\infty})/F, \Zp(i)\big)\big)$ is also finite. 

It therefore remains to show that $H_1\big(H,  H^2_{\Iw}\big(F_\infty/F, \Zp(i)\big)\big)_{\Ga}$ is finite. For this, we require a preliminary lemma.

\bl \label{extension}
Let 
\[ 0\lra M'\lra M\lra M'' \lra 0\]
be a short exact sequence of finitely generated $\Zp\ps{\Gamma}$-modules. Then 
$M_\Gamma$ is finite if and only if $M'_\Gamma$ and $M''_\Gamma$ are finite.

In particular, if $M$ is a finitely generated $\Zp\ps{\Gamma}$-module with $M_\Gamma$ being finite, then $N_\Gamma$ is finite for every subquotient $N$ of $M$.
\el 

\bpf
Indeed, one has the following exact sequence
\[ H_1(\Gamma,M'')\lra M'_\Gamma\lra M_\Gamma \lra M''_\Gamma \lra 0.\]
Clearly, if $M'_\Gamma$ and $M''_\Gamma$ are finite, then so is $M_\Gamma$. Now suppose that $M_\Gamma$ is finite, then so is $M''_\Gamma$. By \cite[Proposition A.1.7]{CS}, the latter implies that $H_1(\Gamma,M'')$ is also finite. Hence it follows from this and the above exact sequence that $M'_\Gamma$ is finite.

The second assertion is immediate from the first.
\epf

We turn back to continue with the proof of our theorem. Applying Proposition \ref{Kato}, we see that there is a
$\Zp\ps{\Gamma}$-subquotient $T$ of
\[ H^2_{\Iw}\big(F_\infty/F, \Zp(i)\big)_H\cong H^2_{\Iw}\big(F(\mu_{p^\infty})/F, \Zp(i)\big)\]
such that $H_1\big(H,  H^2_{\Iw}\big(F_\infty/F, \Zp(i)\big)\big)$ is a successive extension of $T(n_k)$ for some $n_k\geq 1$ with $1\leq k\leq r$. In view of this and taking the first assertion of Lemma \ref{extension} into account, the verification of the finiteness of $H_1\big(H,  H^2_{\Iw}\big(F_\infty/F, \Zp(i)\big)\big)_{\Ga}$ is reduced to showing that $T(k)_\Gamma$ is finite for every $k\geq 1$. Via the second assertion of Lemma \ref{extension}, the latter is reduced to showing that 
\[ \Big(H^2_{\Iw}\big(F_\infty/F, \Zp(i)\big)_H (k)\Big)_{\Gamma}\]
is finite for every $k\geq 1$. Now observe that
\[ \Big(H^2_{\Iw}\big(F_\infty/F, \Zp(i)\big)_H (k) \Big)_{\Gamma}\cong \Big(H^2_{\Iw}\big(F(\mu_{p^\infty})/F, \Zp(i)\big)(k)\Big)_{\Gamma} \cong \Big(H^2_{\Iw}\big(F(\mu_{p^\infty})/F, \Zp(i+k)\big)\Big)_{\Gamma}, \]
where we note that the second isomorphism follows from an application of \cite[Lemma 2.5.1(c)]{Sh22}. By virtue of Propositions \ref{K2 = H2} and \ref{gal descent}, the latter is precisely
\[ H^2\big(G_S(F), \Zp(i+k)\big) \cong K_{2(i+k)-2}(\Op_{F,S}) \]
which is finite as required. This completes the proof of Theorem \ref{Klimitthm} under assumption (b) which also concludes the paper.

\subsection*{Acknowledgments}
The author likes to thank the referee for the many helpful comments
and suggestions on the manuscript.


\footnotesize
\begin{thebibliography}{00}

\bibitem{Bo} A. Borel, Stable real cohomology of arithmetic groups, Ann. Sci. \'Ecole Norm. Sup. (4) 7 (1974), 235--272.

\bibitem{CS} J. Coates and R. Sujatha, Galois cohomology of elliptic curves, Tata Institute of Fundamental
Research Lectures on Mathematics 88.

\bibitem{DSMS} J. Dixon, M. P. F. Du Sautoy, A. Mann and D. Segal, Analytic Pro-$p$ Groups,
2nd edn, Cambridge Stud. Adv. Math. 38, Cambridge Univ.
Press, Cambridge, UK, 1999.

\bibitem{FK} T. Fukaya; K. Kato, A formulation of conjectures on $p$-adic zeta functions in noncommutative Iwasawa theory, Proceedings of the St. Petersburg Mathematical Society, Vol. XII, 1-85, Amer. Math. Soc. Transl. Ser. 2 219, Amer. Math. Soc., Providence, 2006.

\bibitem{GW} K. R. Goodearl and R. B. Warfield, An
introduction to non-commutative Noetherian rings, London Math. Soc.
Stud. Texts 61, Cambridge University Press, 2004.

\bibitem{Har} M. Harris, Correction to $p$-adic representations arising from
descent on abelian varieties, Compos. Math. 121 (2000) 105-108.

\bibitem{Iw73} K. Iwasawa, On ${\bf Z}\sb{l}$-extensions of algebraic number fields, Ann. of Math. (2) 98 (1973), 246-326.

\bibitem{Jannsen89} U. Jannsen, Iwasawa modules up to isomorphism. Algebraic number theory, 171-207, Adv. Stud. Pure Math., 17, Academic Press, Boston, MA, 1989.

\bibitem{Ka2} K. Kato, Universal norms of $p$-units in some
non-commutative Galois extensions, John H. Coates' Sixtieth
Birthday. Doc. Math. 2006, Extra Vol., 551-565.

\bibitem{Kol} M. Kolster, $K$-theory and arithmetic. Contemporary developments in algebraic $K$-theory, 191--258, ICTP Lect. Notes, XV, Abdus Salam Int. Cent. Theoret. Phys., Trieste, 2004.

\bibitem{LT} K. F. Lai and K.-S. Tan,  A generalized Iwasawa's theorem and its application, Res. Math. Sci. 8 (2021), no. 2, Paper No. 20, 18 pp.

\bibitem{Lam} T. Y. Lam, Lectures on Modules and Rings.
Grad. Texts in Math. 189, Springer, 1999.



\bibitem{LimKgroups} M. F. Lim,  On the growth of even $K$-groups of rings of integers in $p$-adic Lie extensions, Israel J. Math. 249 (2022), no. 2, 735-767.

\bibitem{LimKlimit} M. F. Lim, Comparing direct limit and inverse limit of even $K$-groups in multiple $\Zp$-extensions, J. Th\'eor. Nombres Bordeaux 36 (2024) no. 2, pp. 537-555.

\bibitem{LimSh} M. F. Lim and R. Sharifi,
Nekov\'a\v{r} duality over $p$-adic Lie extensions of global fields, Doc. Math. 18 (2013), 621--678.

\bibitem{Ne} J. Nekov\'a\v{r}, Selmer complexes, Ast\'erisque No. 310 (2006), viii+559 pp.

\bibitem{NSW} J. Neukirch, A. Schmidt and K. Wingberg,
Cohomology of Number Fields. 2nd edn., Grundlehren Math.
Wiss. 323 (Springer-Verlag, Berlin, 2008).

\bibitem{Neu} A. Neumann, Completed group algebras without zero divisors,
Arch. Math. 51(6) (1988) 496-499.

\bibitem{PR} B. Perrin-Riou, Arithm\'etique des courbes elliptiques et théorie d'Iwasawa, M\'em. Soc. Math. France (N.S.) No. 17 (1984), 130 pp.


\bibitem{Sh22}    R. Sharifi, Reciprocity maps with restricted ramification, Trans. Amer. Math. Soc. 375 (2022), no. 8, 5361-5392.

\bibitem{Sou} C. Soul\'e, $K$-th\'eorie des anneaux d'entiers de corps de nombres et cohomologie \'etale, Invent. Math. 55 (1979), no. 3, 251--295.

\bibitem{Qui73a} D. Quillen, Higher algebraic $K$-theory. I. Algebraic $K$-theory, I: Higher $K$-theories (Proc. Conf., Battelle Memorial Inst., Seattle, Wash., 1972), pp. 85--147. Lecture Notes in Math., Vol. 341, Springer, Berlin 1973.

\bibitem{Qui73b} D. Quillen, Finite generation of the groups $K\sb{i}$ of rings of algebraic integers. Algebraic $K$-theory, I: Higher $K$-theories (Proc. Conf., Battelle Memorial Inst., Seattle, Wash., 1972), pp. 179--198. Lecture Notes in Math., Vol. 341, Springer, Berlin, 1973.

\bibitem{Vau} D. Vauclair,  Sur la dualit\'e et la descente d'Iwasawa, Ann. Inst. Fourier Grenob. 59 (2009), no. 2, 691-767.
    
\bibitem{V02} O. Venjakob, On the structure theory of the Iwasawa algebra
of a $p$-adic Lie group, J. Eur. Math. Soc. 4(3)
(2002) 271-311.

\bibitem{Vo} V. Voevodsky, On motivic cohomology with $\mathbf Z/l$-coefficients, Ann. of Math. (2) 174 (2011), no. 1, 401--438.

\bibitem{Wei09} C. Weibel, The norm residue isomorphism theorem, J. Topol. 2 (2009), no. 2, 346--372.

\bibitem{WeiKbook} C. Weibel, The $K$-book. An introduction to algebraic $K$-theory. Graduate Studies in Mathematics, 145. American Mathematical Society, Providence, RI, 2013. xii+618 pp.

\end{thebibliography}
\end{document}